%%%%%%%%%%%%%%%%%%%%%%%%%%%%%%%%%%%%%%%%%%%%%%%%%%%%%%
% Some adjunction properties of ample vector bundles % 
%                                                    %
%   by Hironobu Ishihara                             %
%%%%%%%%%%%%%%%%%%%%%%%%%%%%%%%%%%%%%%%%%%%%%%%%%%%%%%

\input amstex
\documentstyle{amsppt}

\vcorrection{-0.5cm}
\hcorrection{+0.3cm}
\pageheight{20.7cm}
\pagewidth{13.0cm}
\magnification=\magstep1
\parskip=6pt

\define\ce{{\Cal E}}
\define\cf{{\Cal F}}

\define\co{{\Cal O}}

\define\pp{{\Bbb P}}
\define\qq{{\Bbb Q}}

\define\cc{{\Bbb C}}
\define\zz{{\Bbb Z}}

\define\pe{{\Bbb P(\Cal E)}}

\define\exc{\operatorname{Exc}}

\define\pic{\operatorname{Pic}}
\define\rk{\operatorname{rank}}
\define\sing{\operatorname{Sing}}

\topmatter

\title 
  Some adjunction properties of ample vector bundles
\endtitle

\author
  Hironobu Ishihara$^{*}$
\endauthor

\leftheadtext{Hironobu Ishihara}
\rightheadtext{Ample vector bundles}

\affil 
  Department of Mathematics,         
  Tokyo Institute of Technology,   
  Oh-okayama, Meguro,    
  Tokyo 152-8551, Japan \\               
  e-mail: ishihara\@math.titech.ac.jp
\endaffil



\thanks
{$^{*}$Research Fellow of the Japan Society for the Promotion of Science}
\endthanks

\abstract
  Let $\ce$ be an ample vector bundle of rank $r$  
  on a projective variety $X$ with only 
  log-terminal singularities.
  We consider the nefness of adjoint divisors
  $K_X+(t-r)\det\ce$ when $t\ge\dim X$ and $t>r$.
  As a corollary,
  we classify pairs $(X,\ce)$ with $c_r$-sectional genus zero.
\endabstract

\endtopmatter

\document

\subhead
  Introduction
\endsubhead

Let $X$ be a smooth projective variety and
$K_X$ the canonical bundle of $X$.
For the study of $X$, it is useful to consider adjoint bundles
$K_X+tL$, where $t$ is a positive integer and
$L$ is an ample line bundle on $X$.
We refer to the books \cite{BS} and \cite{F0} for the properties 
of $K_X+tL$; it is powerful when $t$ is close to $\dim X$.

Recently, as a natural generalization of adjoint bundles,
many authors have considered $K_X+\det\ce$, where $\ce$ is
an ample vector bundle on $X$.
(We say that a vector bundle $\ce$ is ample if $\co_{\pe}(1)$
is ample on $\pe$.)
In particular, Ye and Zhang \cite{YZ} have given a classification
for pairs $(X,\ce)$ when $\rk\ce\ge n-1$ and $K_X+\det\ce$ is
not nef.
Many other results on $K_X+\det\ce$ are obtained when $\rk\ce$ is
close to $\dim X$.
It seems to be difficult to study the nefness of $K_X+\det\ce$
when $\rk\ce$ is small as compared with $\dim X$.

To overcome this difficulty, in the present paper,
we consider the nefness of $K_X+(t-r)\det\ce$ when $t\ge n=\dim X$.
We mainly use vanishing theorems and an estimate of 
the length of extremal rays, hence our argument works on 
projective varieties $X$ with at worst log-terminal singularities.
Our main result is Theorem 2.5 in which we show that
$K_X+(n-r)\det\ce$ is nef unless $(X,\ce)\cong(\pp^4,\co(1)^{\oplus2})$
when $1<r<n-1$.

As a corollary, we see that the $c_r$-sectional genus of the pairs
$(X,\ce)$ is non-negative and we obtain a classification of $(X,\ce)$
with $c_r$-sectional genus zero in case $X$ is log-terminal.
We note that $c_r$-sectional genus is introduced in \cite{I} 
and studied in case $X$ is smooth (see also \cite{FuI}).

\subhead
  1. Preliminaries
\endsubhead

We work over the complex number field $\cc$.
Varieties are always irreducible and reduced.
The tensor products of line bundles are denoted additively,
while we use multiplicative notation for intersection products.
The numerical equivalence is denoted by $\equiv$.
We denote by $L^{\oplus n}$ the direct sum of $n$-copies of 
a line bundle $L$.
The restriction $L|_Y$ of $L$ to a variety $Y$
is often written as $L_Y$.
We denote by $\qq^n$ a (possibly singular) hyperquadric in $\pp^{n+1}$.
A polarized variety $(X,L)$ is said to be a scroll over 
a variety $W$ if $(X,L)\cong(\pp_W(\ce),\co_{\pe}(1))$ 
for some vector bundle $\ce$ on $W$.
The number $\Delta(X,L):=\dim X+L^{\dim X}-h^0(X,L)$ is called 
the $\Delta$-genus of a polarized variety $(X,L)$.

The following facts are main tools of our argument. 

\proclaim{Proposition 1.1} {\rm (\cite{K, Theorem 1})}
  Let $Y$ be a projective variety with only log-terminal singularities
  and $f:Y\to Z$ a contraction morphism of an extremal ray of $Y$.
  Let $E$ be an irreducible component of 
  $\exc(f):=\{y\in Y|f~\text{is not isomorphic at}~y\}$.
  Then $E$ is covered by a family of rational curves $\{C_i\}$
  such that $f(C_i)$ are points and 
  $-K_Y\cdot C_i\le 2(\dim E-\dim f(E))$.
  Moreover, if $f$ is birational, we have 
  $-K_Y\cdot C_i<2(\dim E-\dim f(E))$.
\endproclaim

\proclaim{Proposition 1.2} 
  {\rm (\cite{Z1, Lemma 1}; see also \cite{Z2, Lemma 1}.)}
  Let $Y$ be as in {\rm (1.1)} and 
  $f:Y\to Z$ a birational contraction morphism of an extremal ray $R$.
  Let $F$ be an irreducible component of some positive-dimensional
  fiber of $f$.
  By taking a desingularization $\varphi:V\to F$ of $F$, we get
  $H^q(V,\varphi^*(-H_{F}))=0$ for any $H\in\pic Y$ with 
  $(K_Y+H)R\le 0$ and $q=\dim F$.
\endproclaim

\subhead
  2. Adjunction properties
\endsubhead

Throughout this section, let $X$ be a projective variety 
with at worst log-terminal singularities, $n=\dim X\ge 2$, and 
let $\ce$ be an ample vector bundle of rank $r$ on $X$.

\proclaim{Theorem 2.1}
  When $r\le n+1$, $K_X+(n+2-r)\det\ce$ is always nef.
  Moreover, $K_X+(t-r)\det\ce$ is always nef when 
  $t\ge n+2$ and $r\le t-1$.
\endproclaim

\proclaim{Theorem 2.2}
  When $r\le n$, $K_X+(n+1-r)\det\ce$ is nef unless
  $(X,\ce)\cong(\pp^n,\co(1))$ or $(\pp^n,\co(1)^{\oplus n})$.
\endproclaim

These theorems are proved later;
now we consider the nefness of $K_X+(n-r)\det\ce$ when $r\le n-1$.

\proclaim{Theorem 2.3} {\rm (cf. \cite{F2, Theorem 3.4})}
  When $r=1$, $K_X+(n-1)\ce$ is nef unless $\Delta(X,\ce)=0$ or
  $(X,\ce)$ is a scroll over a smooth curve.
\endproclaim

\demo{Proof}
  The following argument is almost due to Fujita \cite{F2},
  Andreatta and Wi\'sniewski \cite{AW}.
  By the proof of \cite{F2, Theorem 3.4}, we find that
  (2.3) is true except the following case 
  (we set $L:=\ce$ since $r=1$):
  \roster
  \item"($\ast$)" $K_X+(n-1)L$ is not nef and there exists 
                  a birational contraction morphism
                  $f:X\to Z$ such that $(F',L_{F'})\cong(\pp^{n-1},\co(1))$
                  for the normalization $F'$ of an irreducible component
                  $F$ of some fiber of $f$.
  \endroster
  We show that the case ($\ast$) does not occur.
  We consider the structure of $f$ locally in a neighborhood of $F$.
  Since $\dim F=n-1$ and $K_X+(n-1)L$ is not nef,
  the evaluation morphism $f^*f_*L\to L$ is surjective 
  at every point of $F$ by relative spannedness \cite{AW, Theorem 5.1}.
  Hence we have $(F,L_F)\cong(\pp^{n-1},\co(1))$.
  Applying horizontal slicing \cite{AW, Lemma 2.6} repeatedly,
  we get a birational morphism $\varphi:S\to W$
  such that $S$ is a surface with only log-terminal singularities
  and $(K_S+L_S)C<0$ for an irreducible component $C\cong\pp^1$
  of some fiber of $\varphi$.
  Let $\pi:S'\to S$ be a minimal resolution of $S$ and
  let $C'$ be the strict transform of $C$.
  Then $K_{S'}\cdot C'<-1$ and $C'$ deforms in an at least 
  $1$-dimensional family, which derives a contradiction. \qed
\enddemo

\remark{Remark 2.3.1} 
  Polarized varieties $(X,L)$ with $\Delta(X,L)=0$
  have been classified in \cite{F1}.
\endremark
  
\proclaim{Theorem 2.4} {\rm (cf. \cite{Me, Theorem 2})}
  When $r=n-1$, $K_X+\det\ce$ is nef unless
  $(X,\ce)$ is one of the following:
  \roster
  \item"(i)" $(\pp^n,\co(1)^{\oplus(n-1)})$; 
  \item"(ii)" $(\pp^n,\co(1)^{\oplus(n-2)}\oplus\co(2))$;
  \item"(iii)" $(\qq^n,\co(1)^{\oplus(n-1)})$;
  \item"(iv)" $X\cong\pp_C(\cf)$ for a vector bundle $\cf$
                 of rank $n$ on a smooth curve $C$ and 
                 $\ce|_F=\co_{\pp^{n-1}}(1)^{\oplus(n-1)}$ for every fiber 
                 $F\cong\pp^{n-1}$ of the bundle projection $X\to C$;
  \item"(v)" There exists a very ample line bundle $L$ on $X$
             such that $(X,L)$ is a generalized cone on
             $(\pp^2,\co(2))$ or $(\pp^1,\co(e))$ ($e\ge 3$),
             and $\ce=L^{\oplus(n-1)}$. 
  \endroster
\endproclaim

\remark{Remark 2.4.1} 
  The case (v) is overlooked in \cite{Me, Theorem 2},
  but we can recover it.
  We refer to \cite{BS, (1.1.8)} for generalized cones.
\endremark

\proclaim{Theorem 2.5} 
  When $1<r<n-1$, $K_X+(n-r)\det\ce$ is nef unless
  $(X,\ce)\cong(\pp^4,\co(1)^{\oplus2})$.
\endproclaim

\remark{Remark 2.5.1} 
  This theorem is proved by \cite{I} in case $X$ is smooth.
\endremark

\demo{Proof of Theorems 2.1, 2.2 and 2.5}
  Suppose that $t\ge n$ and $r\le t-1$ and $K_X+(t-r)\det\ce$ is not nef.
  When $r=1$, we have $t\ge n$ and $K_X+(t-1)\det\ce$ is not nef.
  Then we are done by \cite{M1, Proposition 2.1} and (2.3).
  When $r=t-1$, we have $r\ge n-1$ and $K_X+\det\ce$ is not nef.
  Then we are done by \cite{Z2, Theorem 1} and (2.4).
  Thus we may suppose that $1<r<t-1$ in the following.
 
  Let $p:\pp_X(\ce)\to X$ be the bundle projection.
  We set $Y:=\pp_X(\ce)$ and denote by $L$ 
  the tautological line bundle of $Y$.
  We can take an extremal ray $R$ of $Y$ such that 
  $p^*(K_X+(t-r)\det\ce)\cdot R<0$ by an argument similar to
  that in \cite{Z1, Claim IV}. 
  Let $f:Y\to Z$ be a contraction morphism of $R$ and
  let $E$ be an irreducible component of $\exc(f)$.
  By (1.1), there exists a rational curve $C\subset E$ 
  belonging to $R$ such that 
  $$ 
    -K_Y\cdot C\le 2(\dim E-\dim f(E))\le 2n
  $$
  since $p|_F:F\to X$ is a finite morphism for every fiber $F$ 
  of $E\to f(E)$.
  On the other hand, we have 
  $$
    \align -K_Y\cdot C
           &=(rL-p^*(K_X+\det\ce))C \\
           &=r\cdot LC-p^*(K_X+(t-r)\det\ce)\cdot C
                      +(t-r-1)(p^*\det\ce)\cdot C \\
           &>r+(t-r-1)r \\
           &=(t-r)r \\
           &\ge 2(t-2),
    \endalign
  $$
  hence $t=n$ or $n+1$, and $LC=1$ or 2.
  If $LC=2$, we see that $t=n$ and $\dim E-\dim f(E)=n$.

  (2.6) Case $LC=1$.

  We have $(K_Y+sL)C<0$ for $s\le t$.
  We use Zhang's idea in \cite{Z1} and \cite{Z2}.
  If $f$ is birational, by (1.2), 
  $H^q(V,\varphi^*(-sL_{F}))=0$ for $s\le t$, 
  where $\varphi:V\to F$ is a desingularization of 
  an irreducible component $F$ of 
  some positive-dimensional fiber of $f$ and $q=\dim F$.
  We get $\chi(V,\varphi^*(-sL_{F}))=0$ for $1\le s\le t$ 
  by Kawamata-Viehweg vanishing theorem.
  Then it follows that $q=n=t$. 
  Let $\mu:W\to F$ be the normalization that factors $\varphi$.
  We get $(W,\mu^*(L_{F}))\cong(\pp^n,\co(1))$ by using
  \cite{F2, Theorem 2.2}.
  Set $\lambda:=(p|_{F})\circ\mu$.
  Then $\lambda:W\to X$ is a finite surjective morphism.
  We can write $\lambda^*(K_X+(n-r)\det\ce)=\co_{\pp^n}(m)$.
  Let $l$ be a line in $W\cong\pp^n$ such that 
  $\lambda(l)\subset X\setminus\sing X$.
  Then we have $m=\lambda^*(K_X+(n-r)\det\ce)\cdot l\in\zz$.
  Set $C':=\mu_*l$ as a 1-cycle.
  We find that 
  $$ 
    \align (K_Y+sL)C'
           &=\mu^*[(s-r)L+p^*(K_X+\det\ce)]_{F}\cdot l \\
           &\le(s-r)+m-(n-r-1)r \\
           &\le 0
    \endalign
  $$ 
  for $s\le n+1$.
  Since $C'\equiv\alpha C$ for some $\alpha>0$, we get 
  $(K_Y+sL)C\le 0$ for $s\le n+1$.
  Then we infer that $\chi(V,\varphi^*(-sL_{F}))=0$ 
  for $1\le s\le n+1$ as before.
  This is a contradiction, thus $f$ is of fiber type.

  Let $F$ be a general fiber of $f$.
  Since $(K_Y+tL)C<0$, we see that $K_F+tL_F$ is not nef.
  Then $t=n$ and $(F,L_F)\cong(\pp^n,\co(1))$ by 
  \cite{M1, Proposition 2.1}.
  Let $U$ be a smooth open subset of $Z$ such that 
  $f^{-1}(z)\cong\pp^n$ for every $z\in U$.
  Set $V:=f^{-1}(U)$. 
  We see that $f|_V:V\to U$ is a smooth morphism.
  It follows that $V$ is smooth and so is $X$.
  Then we obtain that $(X,\ce)\cong(\pp^4,\co(1)^{\oplus2})$ by (2.5.1).

  (2.7) Case $LC=2$.

  We have $(K_Y+sL)C<0$ for $s\le n-1$.
  If $f$ is of fiber type, then $-(K_F+(n-1)L_F)$ is
  ample for a general fiber $F$ of $f$.
  Note that $\dim F=\dim E-\dim f(E)=n$.
  Using Vanishing theorem, we get $\chi(s):=\chi(F,sL_F)=0$ for
  $-(n-1)\le s\le -1$,
  $\chi(0)=h^0(F,\co_F)=1$ and $\chi(1)=h^0(F,L_F)$.
  Then we find that $\Delta(F,L_F)=0$ by Riemann-Roch theorem.
  Hence $(F,L_F)$ is one of the following (\cite{F1}):
  \roster 
  \item"(a)" $(\pp^n,\co(1))$;
  \item"(b)" $(\qq^n,\co(1))$;
  \item"(c)" a scroll over $\pp^1$;
  \item"(d)" a generalized cone over a smooth subvariety $V\subset F$
             with $\Delta(V,L_V)=0$.
  \endroster
  Then there exists a rational curve $l\subset F$
  such that $L_F\cdot l=1$.
  We see that $C\equiv 2l$ and we get 
  $$ 
    2n\ge -K_Y\cdot C>2r(n-r)\ge 4(n-2), 
  $$ 
  a contradiction.
  Thus $f$ is birational.
  Since 
  $$ 
    2n>-K_Y\cdot C>(n-r+1)r\ge 2(n-1), 
  $$ 
  we find that $r=2$ or $(r,n)=(3,5)$.
  If $(r,n)=(3,5)$, then we have $(p^*\det\ce)\cdot C=3$.
  Set $A:=2L-p^*\det\ce$.
  Since $AC=1$, $A$ is an $f$-ample line bundle on $Y$ and
  we have $(K_Y+sA)C<0$ for $s\le 2n-2=8$.
  Then we get a contradiction by using (1.2) as in (2.6).
  Thus we see that $r=2$.
  Since $\dim E-\dim f(E)=n$, there exists an $n$-dimensional
  irreducible component $F$ of some fiber of $f$.
  Since $\dim Y=n+1$ and $K_Y+(n-1)L$ is not nef,
  we infer that $\Delta(F,L_F)=0$
  from the argument in the proof of \cite{A, Theorem 2.1},
  Then we get a contradiction by the same argument that is used
  when $f$ is of fiber type. \qed
\enddemo

\subhead
  3. A corollary on $c_r$-sectional genus
\endsubhead

\definition{Definition 3.1}
  Let $X$ be an $n$-dimensional normal projective variety and
  $\ce$ an ample vector bundle of rank $r<n$ on $X$.
  The {\it $c_r$-sectional genus} $g(X,\ce)$ of a pair $(X,\ce)$ 
  is defined by the formula
  $$ 
    2g(X,\ce)-2:=(K_X+(n-r)c_1(\ce))c_1(\ce)^{n-r-1}c_r(\ce), 
  $$
  where $K_X$ is the canonical divisor of $X$.
\enddefinition

\remark{Remark 3.1.1}
  Let $(X,\ce)$ be as above.
  When $r=1$, $g(X,\ce)$ is called the {\it sectional genus} 
  of a polarized variety $(X,\ce)$.
  We refer to \cite{F0} for the general properties of sectional genus.
  When $r=n-1$, $g(X,\ce)$ is called the {\it curve genus}
  of a generalized polarized variety $(X,\ce)$.
  We refer to \cite{Ba}, \cite{LMS} and \cite{M2}
  for the properties of curve genus in case $X$ is smooth. 
  We have good properties of $g(X,\ce)$ for general $r<n$ 
  in case $X$ is smooth (see \cite{I} and \cite{FuI}). 
\endremark

\proclaim{Lemma 3.2}
  Let $(X,\ce)$ be as in {\rm (3.1)}.
  Then $g(X,\ce)$ is an integer.
\endproclaim

\demo{Proof}
  Let $\pi:X'\to X$ be a desingularization of $X$.
  We get $g(X',\pi^*\ce)\in\zz$ by an argument in \cite{I}.
  We have
  $$ 
    \align 2g(X',\pi^*\ce)-2
     &=(K_{X'}+(n-r)\pi^*c_1(\ce))(\pi^*c_1(\ce))^{n-r-1}\pi^*c_r(\ce)\\
     &=(\pi_*K_{X'}+(n-r)c_1(\ce))c_1(\ce)^{n-r-1}c_r(\ce)\\
     &=2g(X,\ce)-2,
    \endalign
  $$
  hence $g(X,\ce)=g(X',\pi^*\ce)\in\zz$. \qed
\enddemo

As corollaries of Theorems 2.3, 2.4 and 2.5,
we obtain the following theorems.

\proclaim{Theorem 3.3} {\rm (cf.~\cite{F2, Corollary 3.8})}
  Let $L$ be an ample line bundle on a projective variety $X$ 
  with only log-terminal singularities.
  Then $g(X,L)\ge 0$, and $g(X,L)=0$ if and only if 
  $\Delta(X,L)=0$. 
\endproclaim

\demo{Proof}
  First we note that $\Delta(X,L)=0$ implies $g(X,L)=0$
  (see \cite{F1}).
  Assume that $g(X,L)\le 0$. 
  Then $K_X+(n-1)L$ is not nef and it follows that 
  $g(X,L)=\Delta(X,L)=0$ by (2.3). \qed
\enddemo

\proclaim{Theorem 3.4}
  Let $(X,\ce)$ be as in {\rm (3.1)}.  
  Suppose that $2\le r=n-1$ and $X$ has at worst log-terminal
  singularities.  
  Then $g(X,\ce)\ge 0$, and $g(X,\ce)=0$ if and only if
  $(X,\ce)$ is one of the following:
  \roster 
  \item"(i)" $(\pp^n,\co(1)^{\oplus(n-1)})$;
  \item"(ii)" $(\pp^n,\co(1)^{\oplus(n-2)}\oplus\co(2))$;
  \item"(iii)" $(\qq^n,\co(1)^{\oplus(n-1)})$;
  \item"(iv)" $X\cong\pp_{\pp^1}(\cf)$ for a vector bundle $\cf$
                 of rank $n$ on $\pp^1$ and 
                 $\ce|_F=\co_{\pp^{n-1}}(1)^{\oplus(n-1)}$ for every fiber 
                 $F\cong\pp^{n-1}$ of the bundle projection $X\to\pp^1$.
  \item"(v)" There exists a very ample line bundle $L$ on $X$
             such that $(X,L)$ is a generalized cone on
             $(\pp^2,\co(2))$ or $(\pp^1,\co(e))$ ($e\ge 3$),
             and $\ce=L^{\oplus(n-1)}$. 
  \endroster
\endproclaim

\demo{Proof}
  Assume that $g(X,\ce)\le 0$. 
  Then $K_X+\det\ce$ is not nef and $(X,\ce)$ is one of 
  the cases in (2.4).
  In the cases (i), (ii), (iii) and (v) of (2.4), 
  we have $g(X,\ce)=0$. 
  In the case (iv) of (2.4), we have $g(X,\ce)=g(C)$, 
  hence $g(X,\ce)=0$ and $C\cong\pp^1$ by assumption. \qed
\enddemo

\proclaim{Theorem 3.5}
  Let $(X,\ce)$ be as in {\rm (3.1)}.  
  Suppose that $1<r<n-1$ and $X$ has at worst log-terminal
  singularities.  
  Then $g(X,\ce)\ge 0$, and $g(X,\ce)=0$ if and only if
  $(X,\ce)\cong(\pp^4,\co(1)^{\oplus2})$.
\endproclaim

This is shown as in the proof of (3.4) by using (2.5).  

\remark{Acknowledgement}
  The author is grateful to Professors T.~Fujita, M.~Mella
  and Y.~Fukuma for their useful comments.
\endremark

\enddocument